\documentclass[final,a4paper,fleqn,twocolumn]{elsarticle}

\journal{}

\DeclareMathAlphabet{\mathcalold}{OMS}{cmsy}{m}{n}
\DeclareMathAlphabet{\bmathcalold}{OMS}{cmsy}{b}{n}
\DeclareMathAlphabet{\mathbcal}{OMS}{cmsy}{b}{n}

\usepackage{graphicx,mathtools}
\usepackage{multirow} 
\usepackage{subfigure} 
\usepackage{amsmath,amssymb}
\usepackage{hyperref}
\usepackage{color}
\usepackage{tcolorbox}
\tcbuselibrary{xparse}
\usepackage{relsize}
\usepackage{multirow} 
\usepackage{latexsym}
\usepackage{natbib,cuted}
\usepackage{cancel}
\usepackage{array}
\hypersetup{
    colorlinks=true,
    linkcolor=blue,
    filecolor=magenta,      
    urlcolor=cyan,
}
\allowdisplaybreaks

\usepackage{geometry}
\newcommand{\rmd}{\mathrm{d}}
\newcommand{\pp}[2]{\fracd{\partial {#1}}{\partial {#2}}}
\newcommand{\dd}[2]{\fracd{\rmd {#1}}{\rmd {#2}}}
\newcommand{\fracd}[2]{\displaystyle
{\frac{{\displaystyle{#1}}}{{\displaystyle{#2}}}}}

\newcommand{\Qe}{Q_{\mathrm{e}}}

\newcommand{\Qf}{Q_{\mathrm{f}}}

\newcommand{\Qu}{Q_{\mathrm{u}}}

\newcommand{\Xc}{X_{\mathrm{c}}}

\newcommand{\bC}{\boldsymbol{C}}

\newcommand{\bPhi}{\boldsymbol{\Phi}}

\newcommand{\bR}{\boldsymbol{R}}

\newcommand{\bSf}{\boldsymbol{S}_{\mathrm{f}}}

\newcommand{\bS}{\boldsymbol{S}}

\newcommand{\bp}{\boldsymbol{p}}
\newcommand{\bzero}{\boldsymbol{0}}

\newcommand{\dcomp}{d_{\mathrm{comp}}}

\newcommand{\vhs}{v_{\mathrm{hs}}}

\definecolor{ros}{RGB}{148,35,9}   
\usepackage{upgreek}

\newcommand{\bCf}{\bC_{\rm f}}

\newcommand{\br}{\boldsymbol{r}}
\newcommand{\bsigmaC}{\boldsymbol{\sigma}_\mathrm{\!\bC}}
\newcommand{\bsigmaS}{\boldsymbol{\sigma}_\mathrm{\!\bS}}

\def\tsc#1{\csdef{#1}{\textsc{\lowercase{#1}}\xspace}}
\tsc{WGM}
\tsc{QE}

\usepackage{cuted} 

\usepackage[font=small,labelfont=bf,margin=12pt]{caption}
\usepackage{multirow}

\makeatletter
\def\ps@pprintTitle{%
  \let\@oddhead\@empty
  \let\@evenhead\@empty
  \let\@oddfoot\@empty
  \let\@evenfoot\@oddfoot
}
\makeatother

\begin{document}

\begin{frontmatter}

\title{A model of reactive settling of activated sludge: comparison 
 with experimental data}  
 


  \author[1]{Raimund B\"urger}
  \ead{rburger@ing-mat.udec.cl} 
  \author[1]{Julio Careaga}
  \ead{juliocareaga@udec.cl} 
  \author[2]{Stefan Diehl}
  \ead{stefan.diehl@math.lth.se} 
  \author[1]{Romel Pineda\corref{mycorrespondingauthor}}
  \ead{rpineda@ci2ma.udec.cl} 
\cortext[mycorrespondingauthor]{Corresponding author, \texttt{rpineda@ci2ma.udec.cl}}
  

  \address[1]{CI${}^{\mathrm{2}}$MA and Departamento de Ingenier\'{\i}a Matem\'{a}tica, Facultad de Ciencias F\'{i}sicas y Matem\'{a}ticas, Universidad~de~Concepci\'{o}n, Casilla 160-C, Concepci\'{o}n, Chile}
  \address[2]{Centre for Mathematical Sciences, Lund University, P.O.\ Box 118, S-221 00 Lund, Sweden}
  


 \begin{abstract}
A non-negligible part of the biological reactions in the activated sludge process for treatment of wastewater takes place in secondary settling tanks that follow biological reactors.
 It is therefore of interest to develop models of so-called reactive settling that describe 
  the spatial variability of reaction rates caused by the variation of  local concentration 
   of biomass due to hindered settling and compression.  
A reactive-settling model described by a system of nonlinear partial differential equations and a numerical scheme are introduced for the simulation of hindered settling of flocculated particles, compression at high concentrations, dispersion of the flocculated particles in the suspension, dispersion of the dissolved substrates in the fluid, and the mixing that occurs near the feed inlet. 
The model is fitted to experiments from a pilot plant where the sedimentation tank has a varying cross-sectional area.
For the reactions, a modified version of the activated sludge model no.~1 (ASM1) is used with standard coefficients.
The constitutive functions for hindered settling and compression are adjusted to a series of conventional batch settling experiments after the initial induction period of turbulence and reflocculation has been transformed away.
Further (but not substantial) improvements of prediction of experimental steady-state scenarios
  can be achieved   by also fitting additional terms modelling hydrodynamic dispersion.
 \end{abstract}


\end{frontmatter}


\section{Introduction}

\subsection{Scope}

In a water resource recovery facility (WRRF)
wastewater   is mainly  treated through  the activated sludge process (ASP) 
 within a circuit of  biological reactors  coupled with secondary settling tanks (SSTs). 
  The ASP is broadly described in well-known handbooks and monographs 
   \citep{metcalf,droste,chen2020book,makinia20}.  To address  
    some of the general experiences made with  the ASP in real  applications, 
      we mention  that SSTs  contain a substantial amount of the activated sludge and reactions may take place even when oxygen is consumed. In fact, 
up to about one third of the total denitrification has been observed to take  place within  the SSTs \citep{Siegrist1995,Koch1999}. 
An excessive production of nitrogen in an SST, however,  leads to bubbles that destroy the sedimentation properties of the flocculated sludge. On the other hand, 
some denitrification in the SSTs may be preferable, since one can then reduce the nitrate recirculation within the reactors and save pumping energy costs.  
The need to  predict, quantify, 
   and eventually control these and other effects clearly call for the 
    development of   models of so-called reactive settling that include  
  the spatial variability of reaction rates caused by the variation of  local concentration 
   of biomass due to hindered settling and compression. Such a model should depend on time 
    (to handle the transient dynamics of biokinetic reactions within the ASP or the simpler 
     process of denitrification) as well as include some spatial resolution, for instance in
      one space dimension aligned with gravity. From a mathematical point of view such a
       model is naturally posed  by  nonlinear partial differential equations (PDEs).

It is the purpose of this work to introduce a minor extension of the 
  model of reactive settling formulated by nonlinear PDEs 
    by \cite{SDIMA_MOL} to include hydrodynamic dispersion, to present 
      a numerical scheme  for simulations, and first and foremost to 
  calibrate the model  to real data from a pilot WRRF \citep{Kirim2019} where activated sludge reacts with dissolved substrates in an SST whose cross-sectional  varies with depth.
For the biochemical reactions, the activated sludge model no.~1 (ASM1) is used \citep{Henze1987}. 
 This model is slightly adjusted to ensure that only   non-negative concentrations are delivered.
  In fact, 
 the original ASM1  allows  for consumption of ammonia/ammonium when the concentration is zero, 
  which causes unphysical  negative concentrations.
To avoid this, the reaction term for that variable is multiplied by a Monod factor which is close to one for most concentrations and tends to zero fast as the concentration tends to zero.
Standard parameter values for the ASM1 are otherwise used.
For the settling-compression phenomenon, we use a three-parameter 
 constitutive hindered-settling function   and a two-parameter compression function.
Hydrodynamic dispersion of particles and dissolved substrates are each included with one term and a longitudinal dispersivity parameter.
The mixing that occurs near the feed inlet is modelled by a heuristic diffusion term that depends on the volumetric flows through the tank.

As is commonly seen, column batch-settling test data exhibit an initial  so-called 
 induction period   when the (average) settling velocity increases slowly from zero due to initial turbulence and possibly other phenomena.
The induction periods of the data were transformed away with the method by \cite{SDAPNUM1}.
Then the five settling-compression parameters were obtained by a least-squares fit to the transformed batch data.
The obtained batch-settling model is then augmented to include the mixing and the dispersion terms, whose parameters are fitted to one experimental steady-state scenario. Finally, the resulting model is compared to two other steady-state experiments. 

\subsection{Related work}

One-dimensional simulation models based on numerical schemes for the governing PDEs of non-reactive continuous sedimentation in WRRFs have been studied widely in the literature, see e.g.\ \citep{Anderson1981, Chancelier1994, SDsiam2, DeClercq2003, Burger&K&T2005a,DeClercq2008}. 
Non-reactive models addressing the geometry of the vessel as a variable cross-sectional area function include those by \cite{Chancelier1994,SDsiam3,SDcec_varyingA}. 
The relevance of studying the reactions occurring in SSTs has been raised by numerous authors \citep{Hamilton1992,Henze2000ASMbook, Gernaey2006,Kauder2007, Alex2011, Flores2012, Ostace2012, Guerrero2013,Li2013_oxidation,Kirim2019,Freytez2019,Meadows2019}.

We are inspired by the results by  \cite{Kirim2019} and their experimental data, which we gratefully also use in the present article.
Their work extends the so-called B\"urger-Diehl model \citep{SDwatres3,SDCMM2} to the inclusion of the biokinetic ASM1 \citep{Henze2000ASMbook} and a varying cross-sectional area.
Here, we use the model by \cite{SDIMA_MOL} that can be stated as  a system of convection-diffusion-reaction PDEs where the unknowns are the particulate and dissolved components of the ASM1 model.

The simultaneous identification of the hindered-settling and compression functions from experimental data of non-reactive settling, and the removal of the induction period have been studied by \cite{DeClercq2006_thesis, SDAPNUM1}.
More elaborate methods for measuring the concentration of solid particles in the entire tank are presented by \cite{DeClercq2005b, Locatelli2015,Francois2016}.

\section{Materials}

\subsection{Geometry of tank}

\begin{figure*}[t] 
\centering 
\includegraphics[width=0.6\textwidth]{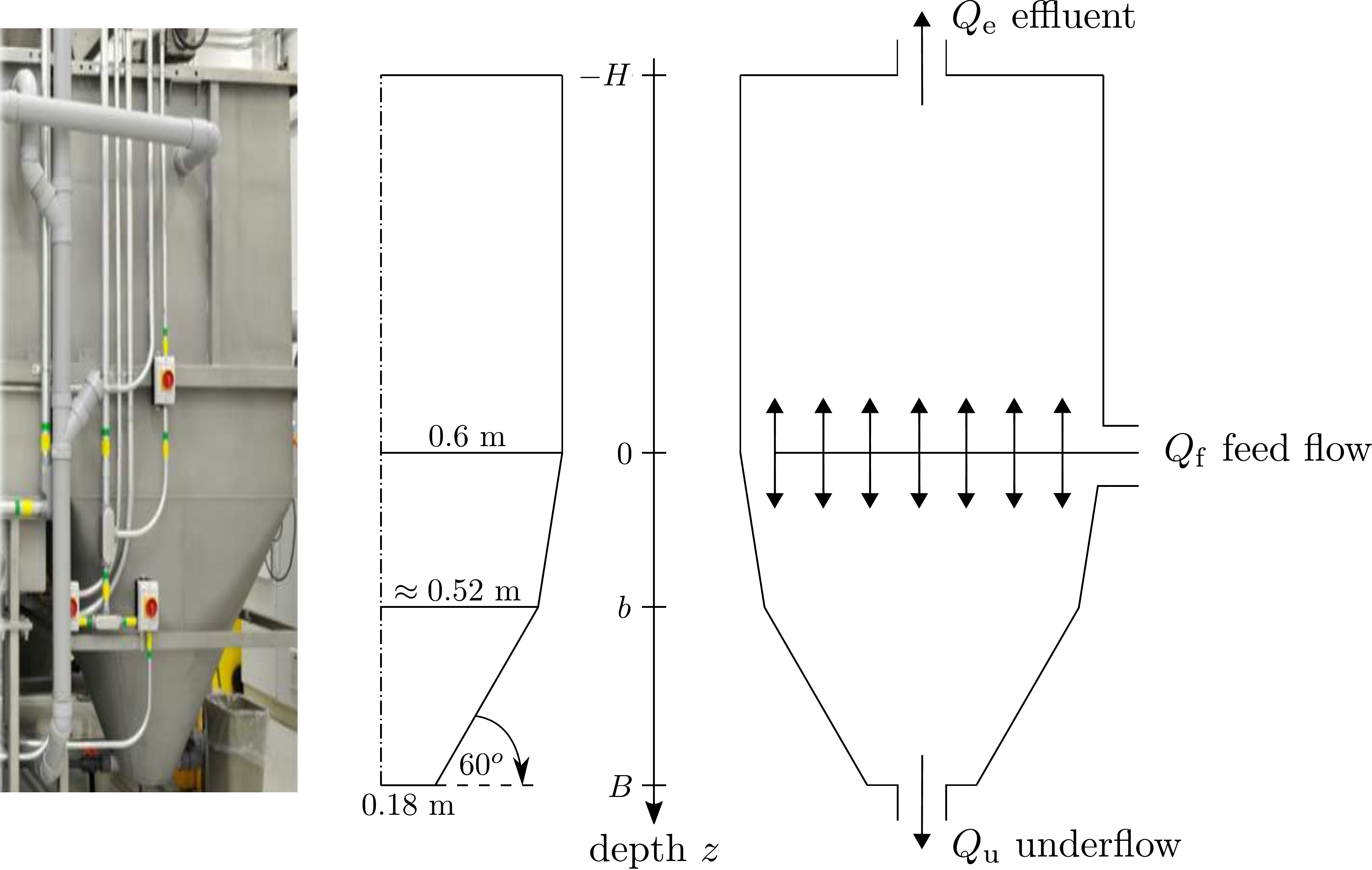}
\caption{Schematic of half of the vertical cross-section of the vessel and the $z$-axis used for the model. The dash-dotted line represents the axis of rotation for the conical part at the bottom.\label{fig:tank}}
\end{figure*}%

A schematic of the tank is shown in Figure~\ref{fig:tank}.
We place a $z$-axis with its origin at the feed level, where the mixture of activated sludge and solubles are fed at the volumetric feed flow $Q_\mathrm{f}$ from the biological reactor.
At the bottom, $z=B=1.1\ \rm m$, there is an  opening through which mixture 
 may leave the unit a controllable 
  volumetric underflow rate~$Q_\mathrm{u}$ and at the top, $z=H=1.25\ \rm m$, the effluent leaves the 
   unit at rate  $Q_\mathrm{e}=Q_\mathrm{f}-Q_\mathrm{u}\geq 0$.
Above the feed inlet, the tank is rectangular with dimensions $1.0\times 1.2~{\rm m}^2=A_0$.
The lowest part, in $b=0.51~{\rm m}<z<B$, is a truncated cone with bottom radius~$r=0.18 \, \mathrm{m}$.
The middle part is has an unknown shape that continuously changes from a (horizontal) rectangle at $z=0$ to a circle with radius $0.52$~m at $z=b$.
 The cross-sectional area in that part is approximated by a convex combination of the rectangle and the circle via 
\begin{align} \label{eq:are_tank} 
A(z)=
	\begin{cases}
	A_0& \text{if $-H \leq z < 0$,} \\
    A_0 + \dfrac{z}{b}\tilde{A}	& \text{if $0 \leq z<b$,} \\[2mm]
	\dfrac{\pi}{3} \bigl( \sqrt{3}r + B-z\bigr)^2	& \text{if $b \leq z< B$,}
	\end{cases}
\end{align}
where 
\begin{align*} 
\tilde{A} =  \dfrac{\pi}{3} \bigl( \sqrt{3}r + B-b\bigr)^2-A_0. 
\end{align*}

\subsection{Activated sludge}

\begin{table*}[t]
\caption{ASM1 variables of the biokinetic reaction model.}
\label{table:AMS1_vari}
\begin{center}  
\begin{tabular}{lll} \hline
Material & Symbol & Unit \\
\hline
Particulate inert organic matter & $X_{\rm I}$  & $\rm (g \ COD)\,m^{-3}$\\
 Slowly biodegradable substrate  & $X_{\rm S}$  & $\rm (g \ COD)\,m^{-3}$\\
 Active heterotrophic biomass	 & $X_{\rm B, H}$& $\rm (g \ COD)\,m^{-3}$\\
 Active autotrophic biomass 	 & $X_{\rm B, A}$& $\rm (g \ COD)\,m^{-3}$\\
 Particulate products  from biomass decay          & $X_{\rm P}$  & $\rm (g \ COD)\,m^{-3}$\\
 Particulate biodegradable   organic nitrogen     & $X_{\rm ND}$ & $\rm (g \ N)\,m^{-3}$\\
 Soluble inert organic matter 				   & $S_{\rm I}$  & $\rm (g \ COD)\,m^{-3}$\\
 Readily biodegradable substrate 				   & $S_{\rm S}$  & $\rm (g \ COD)\,m^{-3}$\\
 Oxygen 										   & $S_{\rm O}$  & $\rm -(g\ COD)\,m^{-3}$\\
 Nitrate and nitrite nitrogen 				   & $S_{\rm NO}$ & $\rm (g \ N)\,m^{-3}$\\
 $\mathrm{NH}_{4}^{+}+\mathrm{NH}_{3}$ nitrogen  & $S_{\rm NH}$ & $\rm (g \ 
 N)\,m^{-3}$\\
 Soluble biodegradable organic nitrogen 		   & $S_{\rm ND}$ & $\rm (g \ N)\,m^{-3}$\\
  Alkalinity 		   & $S_{\rm ALK}$ & $\rm (mol \ CaCO_{3})\,m^{-3}$\\
\hline
\end{tabular} \end{center} 
\end{table*}%

We use the variables of ASM1; see Table~\ref{table:AMS1_vari}, and collect them in the vectors
\begin{align*}
\bC &\coloneqq \bigl(X_\mathrm{I}, X_\mathrm{S}, X_\mathrm{B,H},X_\mathrm{B,A}, X_\mathrm{P},X_\mathrm{ND}\bigr)^{\mathrm{T}},\\ 
\bS &\coloneqq \bigl(S_{\rm I},S_{\rm S}, S_{\rm O},S_{\rm NO},S_{\rm NH},S_{\rm ND}\bigr)^{\mathrm{T}}.
\end{align*}
These concentrations vary with both depth~$z$ from the feed level and the time~$t\geq 0$.
The modified ASM1 is presented in \ref{sec:appendixA} at the end of this paper.

\begin{figure}[t] 
\centering
\includegraphics{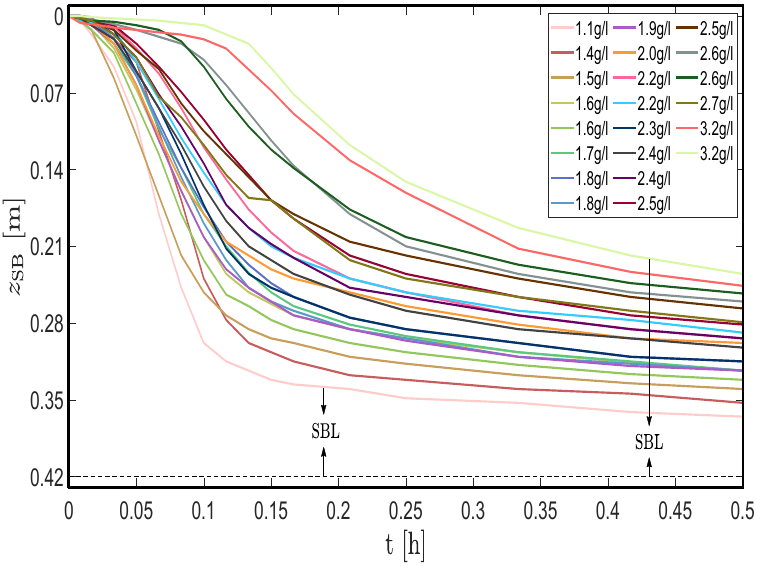}
\caption{Sludge blanket position ($z_{\mathrm{SB}}$) and sludge blanket level (SBL) data obtained from batch settling tests with various initial concentrations~$X_{\mathrm{init}}$.} \label{fig:SBH}
\end{figure}%

\subsection{Batch tests}

The available data contains a series of~22 batch sedimentation tests in cylinders  that were carried out with initial concentrations between $1.1\,\rm g/l$ and $3.2\,\rm g/l$ of activated sludge.
The sludge blanket levels (SBLs), here measured from the bottom, were detected during 30 minutes; see Figure~\ref{fig:SBH}.

\subsection{Steady-state scenarios}

In the measurement campaign of~\cite{Kirim2019}, three operational scenarios were applied to create SBLs at  various  heights: low (L), medium (M) and high (H).
Scenario~L with the lowest SBL, was obtained with the volumetric flows $\Qf=1.0\ {\rm m^3/h}$ and $\Qu=0.5\ {\rm m^3/h}$.
For Scenarios~M and H, the lower values $\Qf=0.65\ {\rm m^3/h}$ and $\Qu=0.15\ {\rm m^3/h}$ were used (lower return flow to the reactor), which create higher SBLs.
In Scenario~H, the internal recycling in the bioreactor is also reduced to obtain a higher nitrate load to the sedimentation tank than in the other scenarios.
The constant in time feed concentrations of the scenarios are shown in Table~\ref{table:CSini&feed}.


\begin{table}[t!]
\caption{Feed concentrations of the solids and substrates of the Scenarios~L, M and H. The units are given in Table~\ref{table:AMS1_vari}.}
\label{table:CSini&feed}
\begin{center}
\begin{tabular}{lcccc} \hline 
& Scenario L& Scenario M & Scenario H \\ 
 \hline 
 $X_{\rm I}$   & 1053.50 &  914.08 & 1309.94 \\
 $X_{\rm S}$   & 46.12 & 40.02 & 57.35\\ 
 $X_{\rm B,H}$  & 1716.60 & 1489.41 & 2134.44 \\ 
 $X_{\rm B,A}$   & 107.71 &  93.45 & 133.92\\ 
 $X_{\rm P}$  &872.55  & 757.08 & 1084.95 \\  
 $X_{\rm ND}$   & 0.04 & 3.30  & 4.73  \\  \hline  
& Scenario L& Scenario M & Scenario H \\ 
\hline 
$S_{\rm I}$  & 14.8 &17.0 & 18.0\\
$S_{\rm S}$   & 0.01  & 0.01  & 0.01 \\ 
$S_{\rm O}$  & 4.5   & 5.2 & 4.48 \\ 
 $S_{\rm NO}$   & 10.95 & 7.0  & 12.65 \\ 
 $S_{\rm NH}$   & 0.022 & 0.01 & 0.021 \\
 $S_{\rm ND}$  &  0.0 & 0.01 & 0.01  \\ \hline  
\end{tabular}
\end{center}
\end{table}%

\section{Methods}

\subsection{Model and numerical method}

The reactive sedimentation model is based on \cite{SDIMA_MOL} with some modifications.
In the derivation of the model, we include hydrodynamic dispersion partly of solubles in the fluid outside the particles, and partly dispersion of the particles in the suspension.
After the governing equations are derived, one may add an ad-hoc term modelling the mixing  near the feed inlet \citep{SDCMM2}.

The solid phase consists of flocculated particles of $k_{\boldsymbol{C}}$ types with mass concentrations~$\smash{C^{(1)},\dots, C^{(k_{\boldsymbol{C}})}}$, which all have the same velocity~$v_X$ and density~$\rho_X$.
The total suspended solids (TSS) concentration is denoted by 
\begin{align}\label{eq:Xtot}
X\coloneqq C^{(1)}+\cdots+ C^{(k_{\boldsymbol{C}})}\quad\text{(mass concentrations)}.
\end{align}
The model is derived in mass concentration units and~\eqref{eq:Xtot} is used during the derivation; however, since the ASM1  is expressed in units that are easier to measure, such as chemical oxygen demand (COD), conversion factors are needed to obtain the true mass concentration.
As is commented on later, the entire model can in fact be used with the units in~Table~\ref{table:AMS1_vari} if the TSS mass concentration is computed by
\begin{multline}\label{eq:X} 
 X \coloneqq \kappa_0 \left(X_\mathrm{I} + X_\mathrm{S} + X_\mathrm{B,H} + X_\mathrm{B,A} + X_\mathrm{P} \right) + X_\mathrm{ND}, \\ 
 \text{where} \quad  \kappa_0=0.75\,\rm g /({\rm g\, COD}).  
\end{multline}
The value of~$\kappa_0$ is taken from \cite[Table~6: WAS, CODp/VSS]{Ahnert2021}.

The liquid phase consists of~$k_{\bS}$ dissolved substrates of concentrations $\smash{S^{(1)} ,\dots, S^{(k_{\bS})}}$ with velocities $\smash{v^{(1)} ,\dots, v^{(k_{\boldsymbol{S}})}}$ and equal density~$\rho_L$.
We collect the unknown concentrations in the vectors~$\bC$ and~$\bS$, and let $\bR_{\bC}$ and $\bR_{\bS}$ denote vectors of the corresponding reaction terms.
Let $\delta(z)$ denote the delta function and $\gamma(z)$ a characteristic function that is one inside the tank and zero outside.
The balance law for each component yields
\begin{multline} \label{eq:balX}  
 A(z)\pp{C^{(k)}}{t} + \pp{}{z}\big(A(z)v_X C^{(k)}\big) \\
 = \delta(z)C_\mathrm{f}^{(k)}(t)\Qf(t) + \gamma(z)A(z)R_{\bC}^{(k)}(\bC,\bS) 
\end{multline} 
for $k=1,\ldots,k_{\bC}$ and 
\begin{multline} \label{eq:balS}  
 A(z)\pp{S^{(k)}}{t} + \pp{}{z}\big( A(z)v^{(k)}S^{(k)} \big) \\
 = \delta(z)S_\mathrm{f}^{(k)}(t)\Qf(t) + \gamma(z)A(z)R_{\bS}^{(k)}(\bC,\bS)
\end{multline} 
for $k=1,\ldots,k_{\bS}$.
It remains to specify the velocities~$v_X$ and~$\smash{v^{(k)}}$ by constitutive assumptions.
As in~\cite{SDIMA_MOL}, one defines the liquid average velocity~$v_L$, the volume average velocity~$q$ of the mixture, and assumes that volume-changing reactions are negligible.
This yields 
\begin{equation}\label{eq:vXvL}
v_X=q+v,\qquad
v_L=q-\frac{X/\rho_X}{1-X/\rho_X}v,
\end{equation}
where the volume average velocity~$q=q(z,t)$ is given via
\begin{align*} 
A(z)q(z,t)=\begin{cases}
-\Qe(t)=\Qu(t)-\Qf(t) & \text{for $z<0$},\\
\Qu(t) & \text{for $z>0$},
\end{cases}
\end{align*}
and the particle excess velocity is given by 
\begin{align*}
v\coloneqq{}& \gamma(z)\left(\vhs(X)-\dcomp(X)\frac{\partial X}{\partial z}+v_\mathrm{disp}\right) \notag\\
={}& \gamma(z)\left(\vhs(X)-\frac{\partial D_{\bC}(X)}{\partial z}+v_\mathrm{disp}\right), 
\end{align*}
where
\begin{equation}\label{eq:dcomp}
D_{\bC}(X)\coloneqq \int_{X_\mathrm{c}}^{X}d_{\rm comp}(\xi)\,\mathrm{d\xi}
\end{equation}
and
\begin{equation*}
d_{\rm comp}(X)\coloneqq
	\begin{cases}
		0 & \text{for $X\leq\Xc$,} \\
		\dfrac{\rho_X\vhs(X)\sigma_\mathrm{e}'(X)}{X g\Delta\rho} & \text{for $X>\Xc$.}
	\end{cases}
\end{equation*}
Here, $\smash{\Delta \rho \coloneqq  \rho_X - \rho_L}$, $g$~is the acceleration of gravity, $\Xc$~is a critical concentration above which the particles are assumed to form a compressible sediment, and $\sigma_\mathrm{e}$ is the  effective solids stress function, which is zero for $X<\Xc$ \citep{Burger&W&C2000}.
The constitutive expressions for~$\sigma_\mathrm{e}$  and the hindered-settling function~$\vhs$ chosen here are defined below.
The velocity~$\smash{v_\mathrm{disp}}$ is due to longitudinal dispersion of particles which is defined via the dispersion flux
\begin{equation}\label{eq:Xdisp}
v_\mathrm{disp}C^{(k)} = -\chi_{\{X<\Xc\}}d_X|q|\pp{C^{(k)}}{z}, \quad k=1,\ldots,k_{\bC},
\end{equation}
where $\chi_{\{X<\Xc\}}$ is a characteristic function, which is one if $X<\Xc$ and zero otherwise, since the particles form a network for higher concentrations.
Furthermore, $d_X$ is the longitudinal dispersivity~[m] of particles in the suspension.
Summing the equalities in~\eqref{eq:Xdisp}, utilizing~\eqref{eq:Xtot} and dividing by~$X$, one gets the formal definition
\begin{multline*} 
 v_\mathrm{disp} =  v_\mathrm{disp}(X,\partial_z X,z,t) \\  \coloneqq  
-\chi_{\{X<\Xc\}}{d}_X|q(z,t)|\pp{\log X}{z}.
\end{multline*}
Each soluble component may also undergo dispersion with respect to the liquid average velocity:
\begin{equation*} 
\big(v^{(k)}-v_L\big)S^{(k)} = -d_L|v_L|\pp{S^{(k)}}{z}, \quad k=1,\ldots,k_{\bS},
\end{equation*}
where $d_L$ is the longitudinal dispersivity [m] of substrates within the liquid.
On top of all these ingredients, we add still another heuristic  term (mathematically, a diffusion term) containing the coefficient
\begin{equation}\label{eq:dmix}
d_\mathrm{mix}(z,\Qu,\Qe)\coloneqq 
\begin{cases}
\mathcal{E}(z,\Qe)
&\text{for $-\alpha_2\Qe<z<0$},\\ 
\mathcal{E}(z,\Qu) 
&\text{for $0<z<\alpha_2\Qu$,} \\
0 &\text{otherwise},
\end{cases}
\end{equation}
where we define 
\begin{align*}
 \mathcal{E}(z,Q) \coloneqq  \alpha_1(\Qu+\Qe)\exp \left(\dfrac{-z^2/(\alpha_2Q)^2}{1-|z|/(\alpha_2Q)}\right).
\end{align*}
The function~$d_\mathrm{mix}$ accounts for the mixing effect due to the feed inlet, where $\alpha_1$ and $\alpha_2$ are parameters.
The larger~$\Qf=\Qu+\Qe$ is, the larger is the effect, and the width above and below the inlet is influenced by $\Qe$ and $\Qu$, respectively.

Substituting \eqref{eq:vXvL}--\eqref{eq:dmix} into the balance laws~\eqref{eq:balX} and \eqref{eq:balS}, one obtains the system of nonlinear PDEs 
\begin{multline}  \label{eq:modelC} 
A(z)\pp{\bC}{t} + \pp{}{z}\big(A(z) \mathcal{V}_{\boldsymbol{C}}\bC\big) 
\\   
= \pp{}{z}\left(A(z)\gamma(z)\Big(\chi_{\{X<\Xc\}}d_X|q| +d_\mathrm{mix}\Big)\pp{\bC}{z}\right)\\
+ \delta(z)\bCf(t)\Qf(t) + \gamma(z)A(z)\bR_{\bC}(\bC,\bS), 
\end{multline} 
\vspace*{-2em}
\begin{multline}  \label{eq:modelS} 
 A(z)\pp{\bS}{t} + \pp{}{z}\big( A(z) \mathcal{V}_{\boldsymbol{S}} \bS\big) \\ 
  =  
 \pp{}{z}\left(A(z)\gamma(z)\Big(d_L|\mathcal{V}_{\boldsymbol{S}}| + d_\mathrm{mix}\Big)\pp{\bS}{z}\right) \\
 + \delta(z)\bSf(t)\Qf(t) + \gamma(z)A(z)\bR_{\bS}(\bC,\bS),
\end{multline}
where the velocity functions are
\begin{align*}
  \mathcal{V}_{\boldsymbol{C}} ={}& \mathcal{V}_{\boldsymbol{C}}(X,z,t)  \\ 
 \coloneqq {}& q+\gamma(z)\left(\vhs(X)-\frac{\partial D_{\bC}(X)}{\partial z}\right),\\
 \mathcal{V}_{\boldsymbol{S}} ={}& \mathcal{V}_{\boldsymbol{S}}(X,\partial_z X,z,t)  \\
\coloneqq{}&  q-\frac{X\gamma(z)}{\rho_X-X}\biggl(\vhs(X)\\
& -\pp{}{z} \Big(D_{\bC}(X) +\chi_{\{0<X<\Xc\}}d_X|q|{\log X}\Big)\biggr).
\end{align*}
The reaction terms have the forms
\begin{align}\label{eq:R} \begin{split} 
\bR_{\bC}(\bC,\bS) & \coloneqq \bsigmaC\br(\bC,\bS), \\
\bR_{\bS}(\bC,\bS) & \coloneqq \bsigmaS\br(\bC,\bS), \end{split} 
\end{align}
where $\bsigmaC$ and $\bsigmaS$ are dimensionless stoichiometric matrices and the vector $\br$ [kg/(sm$^3$)] contains the processes of biokinetic reactions of carbon and nitrogen removal; see \ref{sec:appendixA}, where the modified ASM1 is presented.
That reaction model is expressed in units of Table~\ref{table:AMS1_vari} rather than the mass concentration units in which the present model is derived.
Conversion factors between units are needed.
As we have shown in~\cite{bcdp_part2}, the structure of the governing PDEs~\eqref{eq:modelC}, \eqref{eq:modelS}, where all the nonlinear coefficients depend on~$X$ and $\partial_z X$, and the equations otherwise are linear in~$\bC$ and $\bS$, the model and numerical scheme can be used straightforwardly also with the COD units if only the formula~\eqref{eq:X} is used.

The constitutive functions of hindered settling and sediment compressibility are chosen as  
\begin{align}
\vhs(X)& \coloneqq  \frac{v_0}{1 + (X/ \bar{X})^\eta},\label{eq:vhs}\\
\sigma_\mathrm{e}(X)& \coloneqq 
	\begin{cases}
		0 & \text{for $X\leq\Xc$,} \\
		\alpha(X-\Xc)\vhs(X) & \text{for $X>\Xc$}
	\end{cases}\label{eq:sigmae}
\end{align}
\citep{SDAPNUM1}, where $v_0$, $\bar{X}$, $\eta$ and $\alpha$ are positive parameters.  

Approximate solutions of the system~\eqref{eq:modelC}, \eqref{eq:modelS} are obtained by the numerical method described in \ref{sec:appendixB}, which is an adjusted version the method by~\cite{SDIMA_MOL}.
For the simulations  we use 100~layers for the spatial discretization of the tank.

\subsection{Preparation of batch-test data by removing the induction period}

As in most  batch settling tests presented in the literature, those in Figure~\ref{fig:SBH} show an initial induction period when several phenomena occur, such as turbulence because of mixing before $t=0$ to obtain a homogeneous initial concentration~$X_\mathrm{init}$, reflocculation of broken flocs and possibly rising air bubbles.
Such phenomena are not captured with the assumptions made above.
A PDE model, with the constitutive functions~\eqref{eq:vhs} and \eqref{eq:sigmae}, for the TSS concentration~$X$ during batch sedimentation in a column with constant cross-sectional area is ($z$~is the depth from the top of the column)
\begin{equation}\label{eq:epd}
\pp{X}{t} + \pp{}{z}\big(\vhs(X)X\big)=\pp{}{z}\left ( d_{\rm comp}(X)\pp{X}{z} \right ).
\end{equation}
To fit this to the batch data, we first transform away those phenomena with the technique by \cite{SDAPNUM1} before we calibrate the parameters of the constitutive functions.
The batch tests revealed that for $X >  3.2 \, \mathrm{kg}/\mathrm{m}^2$, the sludge showed a compressible behaviour, wherefore the critical concentration was set to $X_\mathrm{c}=3.2$~kg/m$^2$.
Other parameters are $\rho_X = 1050\, \rm kg/m^3$, $\Delta\rho = 52\, \rm kg/m^3$ and $g = 9.81 \,\rm m/s^2$. 

Let $z=z_\mathrm{p}(t)$ be the path of a solid particle at the SBL that starts at depth $z_\mathrm{p}(0)$.
In a suspension with ideal particles, e.g.\ glass beads, and $X_\mathrm{init}<X_\mathrm{c}$, it is well known that the sludge blanket initially decreases at a constant velocity $v_\mathrm{hs}(X_\mathrm{init})$.
During the induction period, the velocity of each particle increases from zero to its maximum velocity $v_\mathrm{hs}(X_\mathrm{init})$.
Thus, the velocity can be written 
\begin{equation}\label{eq:velpart}
z_\mathrm{p}'(t) = G(t)v_\mathrm{hs}(X_\mathrm{init})
\end{equation}
with a function~$G$ that satisfies $G(0)=0$ and increases to one at the end of the induction period.
For the sludges investigated by \cite{SDAPNUM1}, the  function 
\begin{equation}\label{eq:g}
G(t)= 1 - \exp \bigl( -(t/\bar{t})^p \bigr)
\end{equation}
was appropriate. Here $\bar{t}$~and $p$~are parameters that depend on~$X_{\rm init}$.
Thus, we have to fit such a function $G$ for each batch test.
The particle path is obtained by integrating \eqref{eq:velpart}:
\begin{multline}\label{eq:particle_path} 
z_\mathrm{p}(t) =z_\mathrm{p}(0) \\ 
 \quad + v_\mathrm{hs}(X_{\rm init})\int_{0}^{t}\bigl( 1-\exp\bigl( -(s/\bar{t})^p\bigr) \bigr) \, \mathrm{d}s.
\end{multline}
Under the change of the time coordinate 
\begin{equation}\label{eq:tau}
\tau \coloneqq  \int_{0}^{t} g(s) \, \mathrm{d}s, \quad t>0,
\end{equation}
then the data plotted against~$\tau$ (instead of~$t$) will be nearly straight lines -- the induction period has been transformed away. 
We refer to \cite{SDAPNUM1} for justification.

\begin{figure*}[t]  
\centering 
\includegraphics{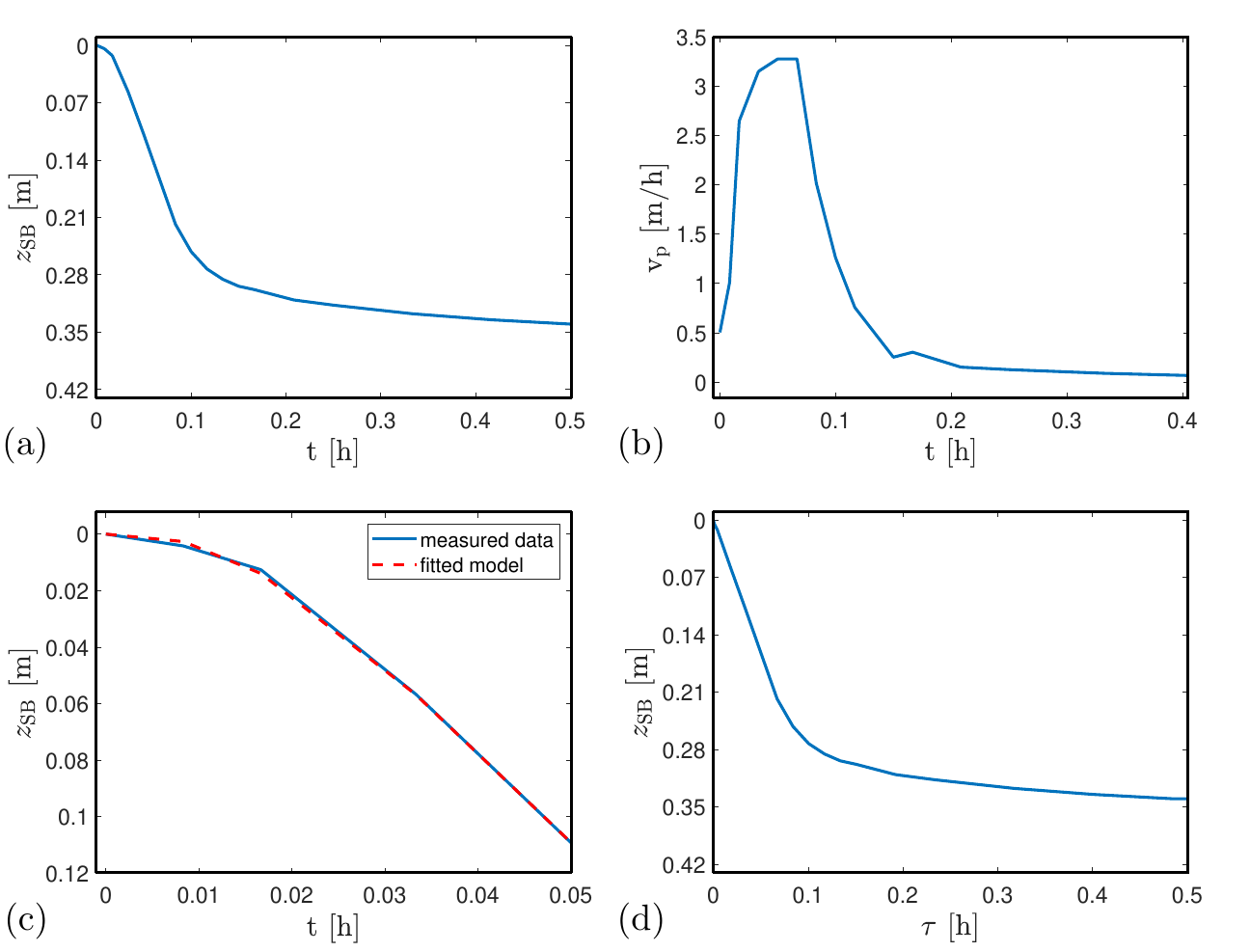} 
\caption{(a) A batch test with $X_\mathrm{init}= 1.5\, \rm kg/m^3$. The curve $\smash{z=z_\mathrm{p}(t)}$ is the location of a particle at the SBL.
(b) Estimated velocity $\smash{v_\mathrm{p}(t):=z_\mathrm{p}'(t)}$ of the curve shown in a.
(c) The initial induction period of the curve in plot a and the fitted model~\eqref{eq:particle_path}. 
(d) The resulting batch sedimentation curve $\smash{z=z^\mathrm{data}(\tau)}$ after rescaling  time with~\eqref{eq:tau}.} \label{fig:path_particle}
\end{figure*}%

Figure~\ref{fig:path_particle}~(a) shows a batch test with an initially concave SBL curve during the induction period.
To determine when that period ends, we take central finite differences of the data to obtain the velocity; see Figure~\ref{fig:path_particle}~(b).
The maximum velocity in that example occurs at approximately $t=0.05\;$h, i.e., we should have $G(0.05\;{\rm h})\approx 1 $.
To data from the induction time interval one performs a nonlinear least-squares fit of the function $z_\mathrm{p}(t)$ to find the parameters $\bar{t}$ and $p$.
The parameters in that example are (see the curve in Figure~\ref{fig:path_particle}~(c))
\begin{align*}
v_\mathrm{hs}(X_{\rm init})
&=8.9261\times10^{-4}\,\mathrm{m}/\mathrm{s},\\
\bar{t}&=64.63\,\mathrm{s}\approx 1.795\times 10^{-2}\,\mathrm{h}, \\
p&=1.641. 
\end{align*}
With these parameters, the rescaling of the time variable by~\eqref{eq:tau} produces the new path of the sedimentation curve in Figure~\ref{fig:path_particle}~(d).

\subsection{Calibration of the settling-compression model}\label{sec:calibvhs}

\begin{figure}[t] 
\centering 
\includegraphics{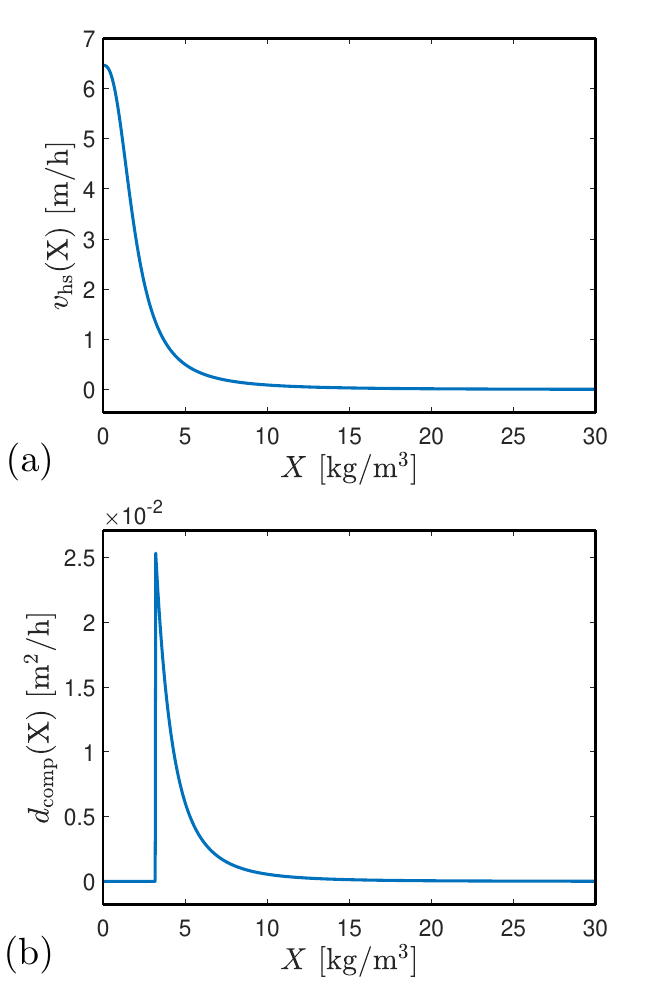}
\caption{Graphs of the constitutive functions: (a) hindered settling velocity, 
 (b)  compression.}
\label{fig:fConstitu}
\end{figure}

\begin{figure}[t!] 
\centering 
\includegraphics{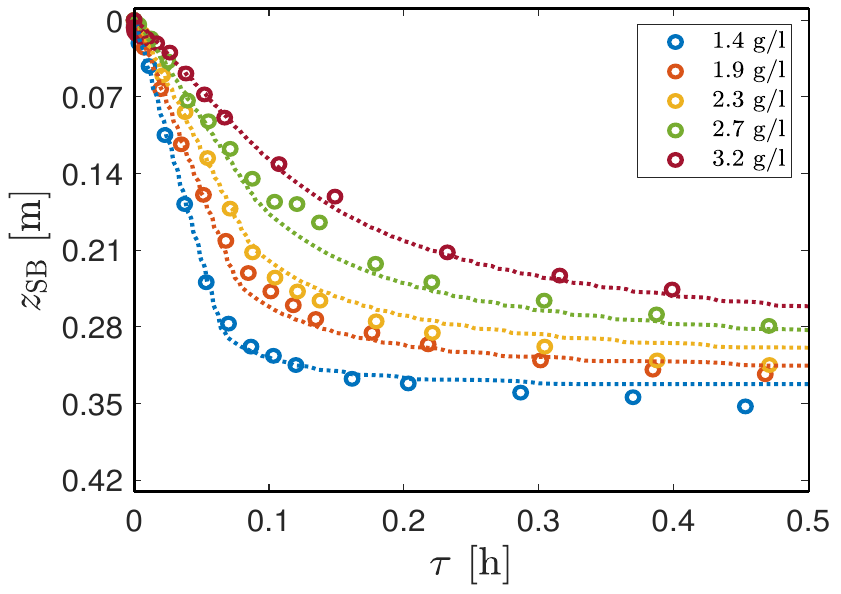}
\caption{Selected 
sludge blanket level (SBL) trajectories 
 in transformed time. Measured data (circles) compared  with 
   the optimized parameters in model \eqref{eq:epd} (dotted curves).} \label{fig:Dat_Sim}
\end{figure}

We denote the transformed trajectories of the batch sedimentation curves by
\begin{equation*}
z=z_j^\mathrm{data}(\tau_i),\quad i=1,\dots,N_j,
\quad j=1,\dots,N_\mathrm{exp},
\end{equation*} 
where $N_j$ is the number of data points in experiment~$j$, and $N_\mathrm{exp}$ is the number of batch experiments.
Now, we proceed to find the optimal parameters ($ v_0 $, $ \bar {X} $, $ \eta $ and $ \alpha $) for the constitutive functions in \eqref{eq:vhs} and \eqref{eq:sigmae}. 
This is done by minimizing the sum of squared errors
\begin{multline}\label{eq:SSE} 
 E(v_0,\bar{X},\eta,\alpha) \\ 
 \coloneqq  \sum_{j=1}^{N_{\rm exp}}\sum_{i=1}^{N_{j}} \left ( {z_j^{\rm data}(\tau_i) - \hat{z}_j(\tau_{i};v_0,\bar{X},\eta,\alpha) } \right)^{2},
\end{multline} 
where $\hat{z}_j(\tau_{i};v_0,\bar{X},\eta,\alpha)$ is the estimated value obtained by numerical simulation of the model~\eqref{eq:epd} with the method by~\cite{SDcec_varyingA} with 100~spatial layers.
The optimal parameters after minimizing the objective function \eqref{eq:SSE} with the robust Nelder-Mead simplex algorithm \citep{nelder1965simplex} are
\begin{alignat*}2
v_0 &= 6.46\, {\rm m/h}, &\qquad& \bar{X} = 1.89\, {\rm kg/m^3},\\
\eta &= 2.55, &\qquad& \alpha = 381605.95\,\rm m^2/h^2.
\end{alignat*}
Graphs of the constitutive functions with these parameter values are shown in Figure~\ref{fig:fConstitu}.
Finally, we show the outcome of the optimization by comparing some of the batch tests with the corresponding simulations in Figure~\ref{fig:Dat_Sim}.

\subsection{Calibration of the dispersion and mixing parameters}\label{sec:calibdisp}

In an attempt to improve the model calibrated from the batch experiments, we now use the experimental steady-state solution of Scenario~M with its four data points of each concentration (see Figure~\ref{fig:Scen2c}) to calibrate the remaining parameters of the model
\begin{align*} 
\bp\coloneqq (d_X, d_L, \alpha_1, \alpha_2)
\end{align*} 
with the error function  
\begin{align} \label{eq:Edisp} \begin{split} 
& E_\mathrm{disp}(\bp)  \coloneqq \sum_{j=1}^{4}
\Biggl( 
\frac{|X^\mathrm{data}(z_j)-\hat{X}(z_j;\bp)|}{\max_i \{X^\mathrm{data}(z_i)\}} \\
&\quad  + \dfrac{|(S^\mathrm{data}_\mathrm{I}+S^\mathrm{data}_\mathrm{S})(z_j) -(\hat{S}_\mathrm{I}+\hat{S}_\mathrm{S})(z_j;\bp)|}{\max_i\{(S^\mathrm{data}_\mathrm{I}+S^\mathrm{data}_\mathrm{S})(z_i)\}}\\
& \quad + \dfrac{|S^\mathrm{data}_\mathrm{O}(z_j)-\hat{S}_\mathrm{O}(z_j;\bp)|}{\max_i \{S^\mathrm{data}_\mathrm{O}(z_i)\}} \\
& \quad + \dfrac{|S^\mathrm{data}_\mathrm{NO}(z_j)-\hat{S}_\mathrm{NO}(z_j;\bp)|}{\max_i \{S^\mathrm{data}_\mathrm{NO}(z_i)\}} \\
& \quad + \dfrac{|S^\mathrm{data}_\mathrm{NH}(z_j)-\hat{S}_\mathrm{NH}(z_j;\bp)|}{\max_i \{S^\mathrm{data}_\mathrm{NH}(z_i)\}}
\Biggr),
\end{split} \end{align}  
where $\hat{X}(z_j;\bp)$ is the steady-state approximation of $X$ of the PDE system~\eqref{eq:modelC}, \eqref{eq:modelS} obtained by the numerical scheme~\eqref{eq:MOL} with 100 spatial layers (see \ref{sec:appendixB}) at the point $z=z_j$, analogously for the other concentrations. 
The variables labelled by ``data'' are the experimental data points from Scenario~M. 

We use the Nelder-Mead algorithm also for the minimization of $E_\mathrm{disp}$, whose function value at~$\bp$ is obtained by simulating to a steady state.
Such simulations are made with the numerical scheme in \ref{sec:appendixB} and need initial conditions preferably close to the expected steady state.
Since the nature of the Nelder-Mead algorithm is to compare only function values at the corners of a simplex and most of the time no large step is taken, we choose the initial conditions in the following way:
\begin{itemize}

\item[i.] Given the starting point~$\bp^0$ of the optimization iteration (e.g.\ $\bp^0 = \boldsymbol{0}$), simulate from any initial data $\smash{\bC_{\rm init}(z)}$ to obtain approximate steady states $\smash{ \boldsymbol{\hat{C}}(z;\bp^0)}$ (analogously for $\bS$), which give the value~$E_\mathrm{disp}(\bp^0)$.

\item[ii.] For optimization iteration $k\geq 1$, i.e., given~$\bp^k$, simulate with the initial data $\smash{\bC_{\rm init}(z) = \boldsymbol{\hat{C}} (z;\bp^{k-1})}$ (analogously for~$\bS$) to obtain~$E_\mathrm{disp}(\bp^k)$.

\end{itemize}
The optimum vector of parameters for Scenario~M is
\begin{equation*}
\bp^\ast=(0.004157,0.03817,0.01678,0.0895),
\end{equation*}
where the units of the parameters $d_X$, $d_L$, $\alpha_1$ and $\alpha_2$ are $[\rm m]$, $[\rm m]$, $[\rm m^{-1}]$ and $[\rm h/m^2]$, respectively.
The corresponding approximate steady-state profiles are shown in Figure~\ref{fig:Scen2c}.
To investigate the impact with and without the mixing term~$d_\mathrm{mix}$, we also carried out the optimization with the reduced vector of parameters 
$\boldsymbol{u}\coloneqq (d_X,d_L)$ (instead of~$\bp$) and found the optimum point
\begin{equation*}
\boldsymbol{u}^{\ast}= (0.07044 , 0.04837)~\mathrm{m}.
 \end{equation*}
The corresponding steady states are also shown in Figure~\ref{fig:Scen2c}.

\section{Results}

\begin{figure*}[t] 
\centering 
\includegraphics{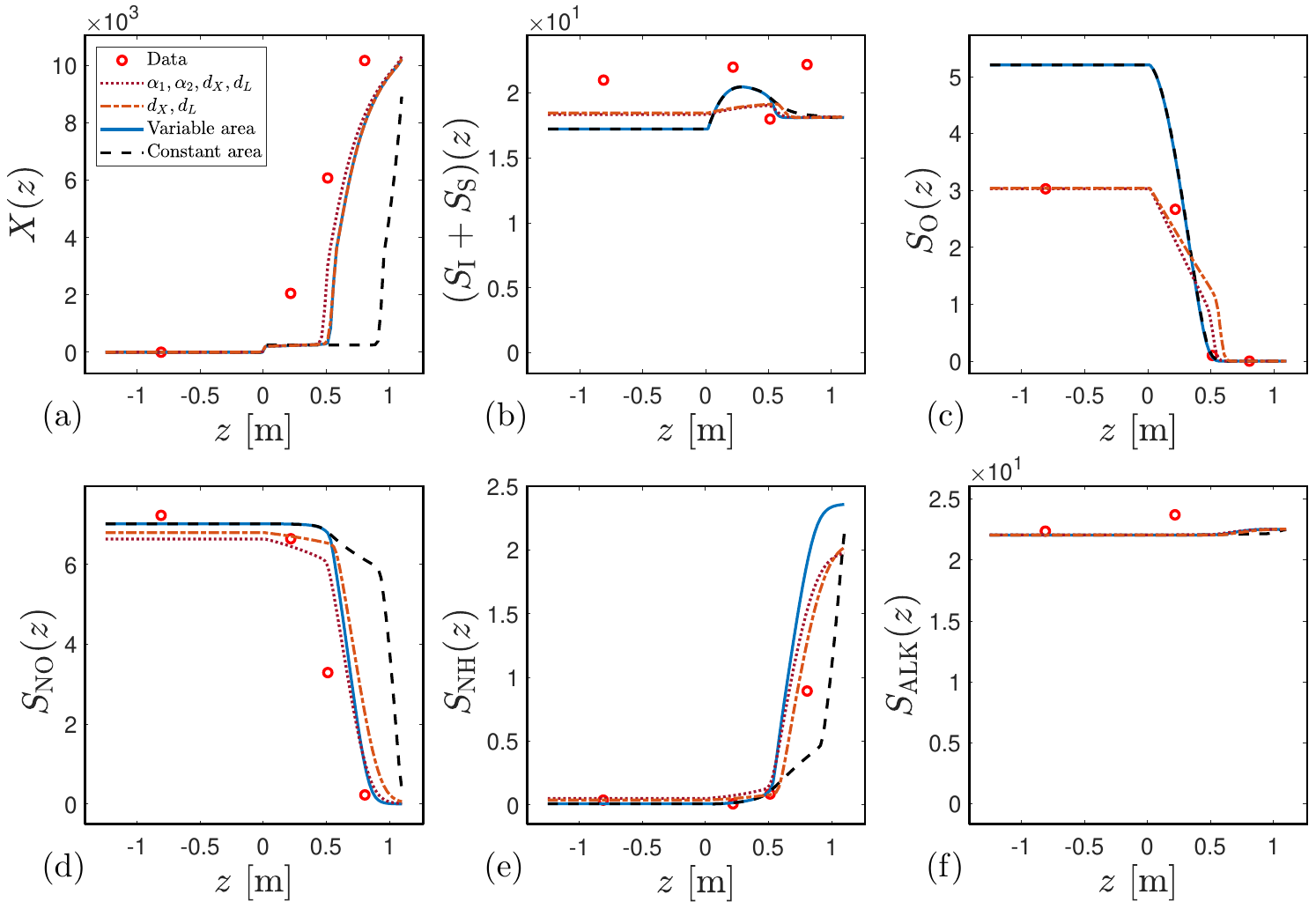}
\caption{Scenario~M: Concentrations with units given in Table~\ref{table:AMS1_vari}: (a) total suspended solids, (b) soluble chemical oxygen demand, (c) oxygen, (d) nitrate and nitrite nitrogen, (e) $\mathrm{NH}_{4}^+  + \mathrm{NH}_3$ nitrogen, (f) alkalinity.
Each plot shows experimental data (circles), simulation results at steady state ($T = 24$~h in Figures~\ref{fig:Scen2a} and \ref{fig:Scen2b}) with
a varying cross-sectional area \eqref{eq:are_tank} (solid line) and constant cross-sectional area (dashed line), both obtained with $\alpha_1=\alpha_2=d_X=d_L=0$,
with optimized parameters $(\alpha_1 , \alpha_2, d_X, d_L )$ (dotted line), and $(d_X, d_L )$ 
(dash-dotted line).}  \label{fig:Scen2c}
\end{figure*}%

Figure~\ref{fig:Scen2c} shows the experimental data of Scenario~M together with simulated steady states of four different models.
The solid blue curves show the simulated steady states with the model calibrated from the batch experiments only; that is, with only the constitutive functions for hindered settling and compression, and neither mixing nor dispersion ($\bp = \boldsymbol{0}$).
As a comparison, the dashed black curves show the results when instead a cylindrical tank is used.
That tank has the same height and volume as the one in Figure~\ref{fig:tank}, which means that the cross-sectional area is $A = 0.9653\,\rm m^2$; otherwise the conditions are the same.
The two remaining curves in each subplot show the steady states when partly only $d_X$ and $d_L$ are fitted (and $\alpha_1=\alpha_2=0$), and partly when all four parameters~$\bp$ has been used in the fitting.

\begin{figure*}[t] 
\centering 
 \includegraphics{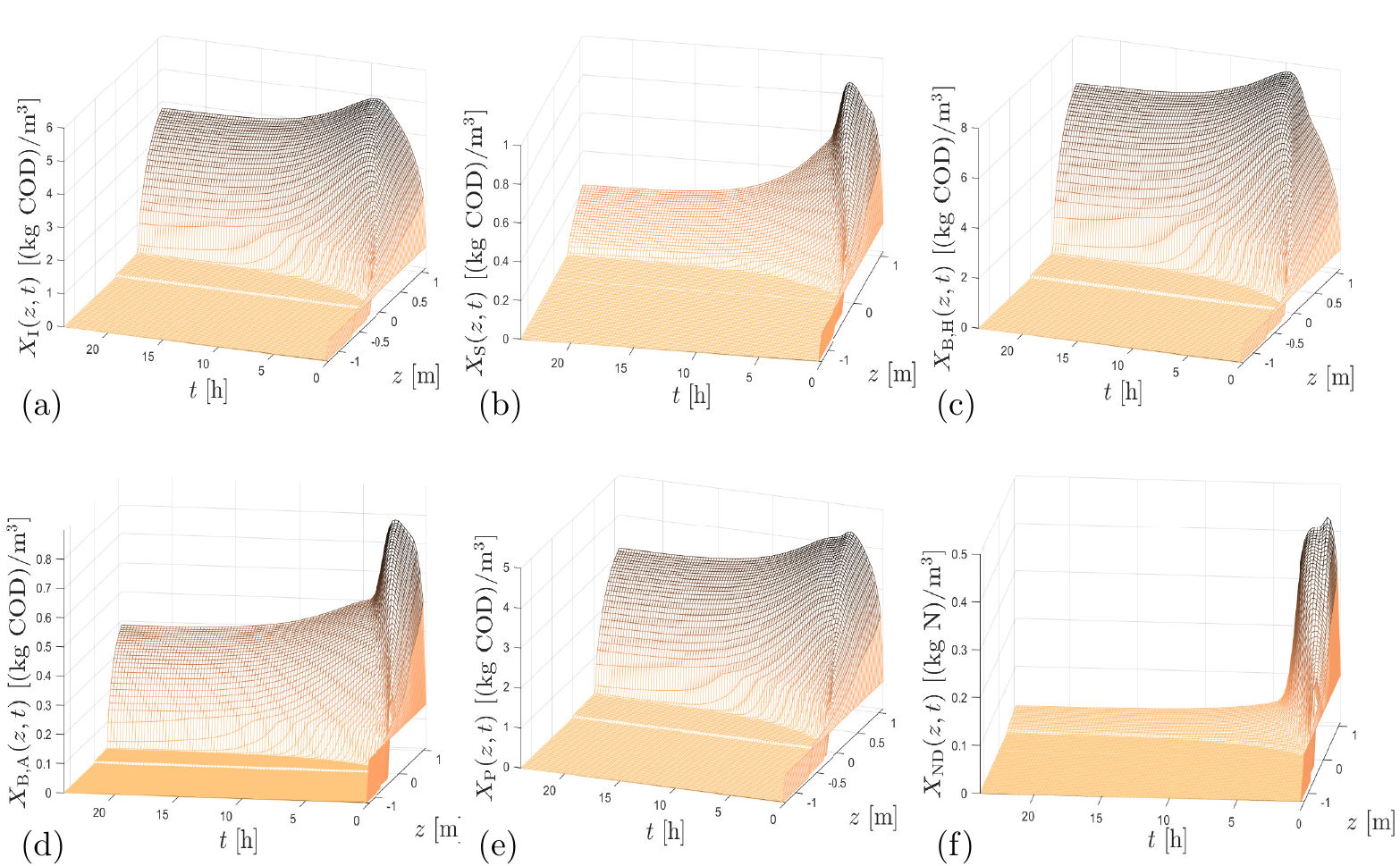} 
 \caption{Scenario~M: Numerical simulation of solids during 24~h:  (a) particulate inert organic matter, 
  (b) slowly biodegradable substrate, (c) active heterotrophic biomass,   (d) 
  active autotrophic biomass,  (e) particle products arising from biomass decay  and 
  (f) particulate biodegradable organic nitrogen.} \label{fig:Scen2a}
\end{figure*}

\begin{figure*}[t] 
\centering 
\includegraphics{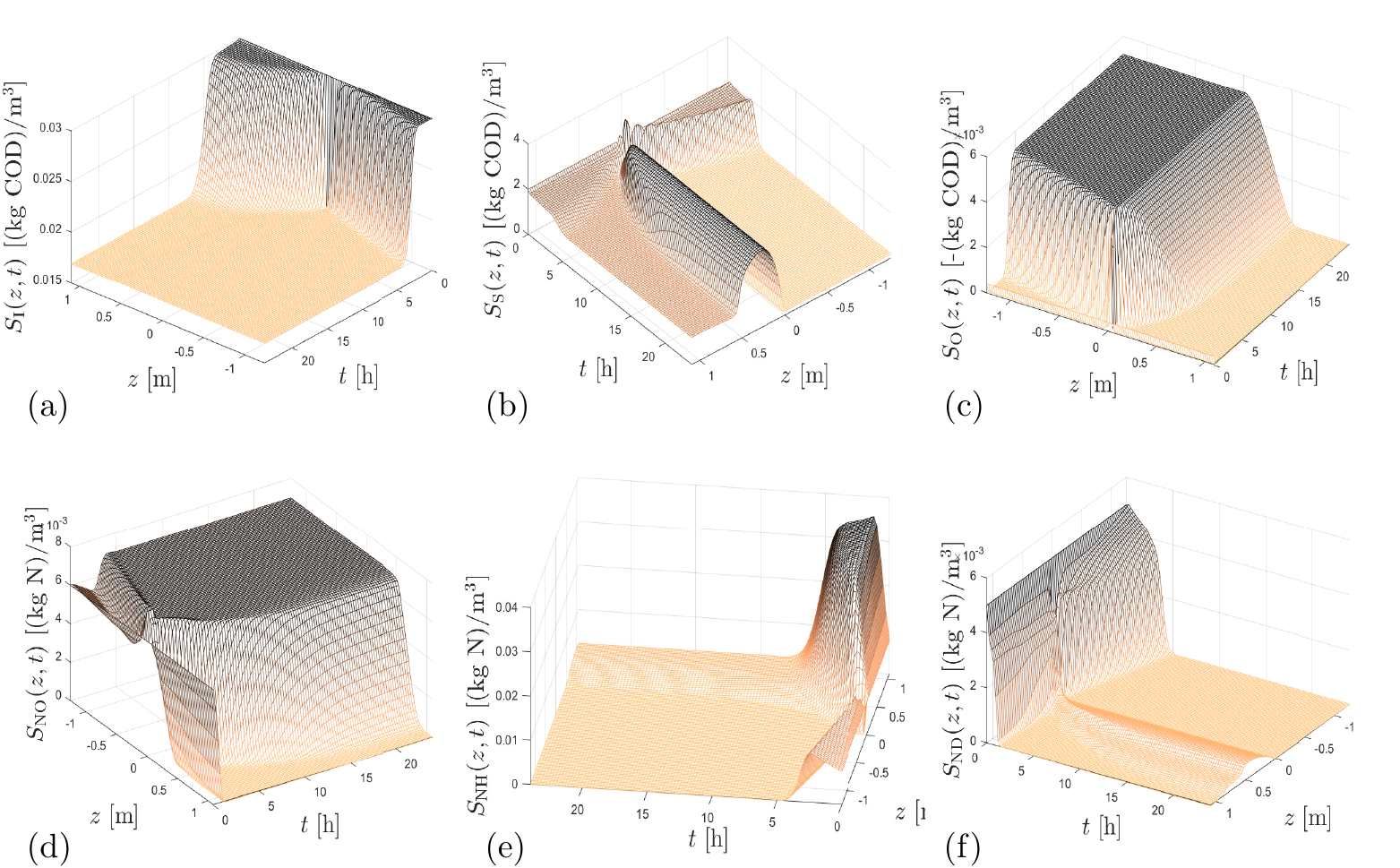}
 \caption{Scenario~M: Numerical simulation of solubles during 24~h: (a) soluble inert organic matter, 
 (b) readily biodegradable substrate, (c) oxygen, (d) nitrate and nitrite nitrogen, 
  (e) $\mathrm{NH}_{4}^+  + \mathrm{NH}_3$ nitrogen, (f) soluble biodegradable organic nitrogen.} \label{fig:Scen2b}
\end{figure*}

For purpose of demonstration, we show in Figures~\ref{fig:Scen2a} and~\ref{fig:Scen2b} a full dynamic simulation of all concentrations (except alkalinity) with the feed input concentrations shown in Table~\ref{table:CSini&feed} and with the constant initial concentrations
\begin{align*}
\bC_\mathrm{init} &= (650,150,800,150,700,100)^\mathrm{T},\\
\bS_\mathrm{init} &= (30.0,2.0,0.4,6.07.5,5.0)^\mathrm{T},
\end{align*}
with units as in Table~\ref{table:AMS1_vari}.
The simulation is performed until~$T=24\,$h, where the solution is in approximate steady state.

\begin{figure*}[t] 
\centering 
 \includegraphics{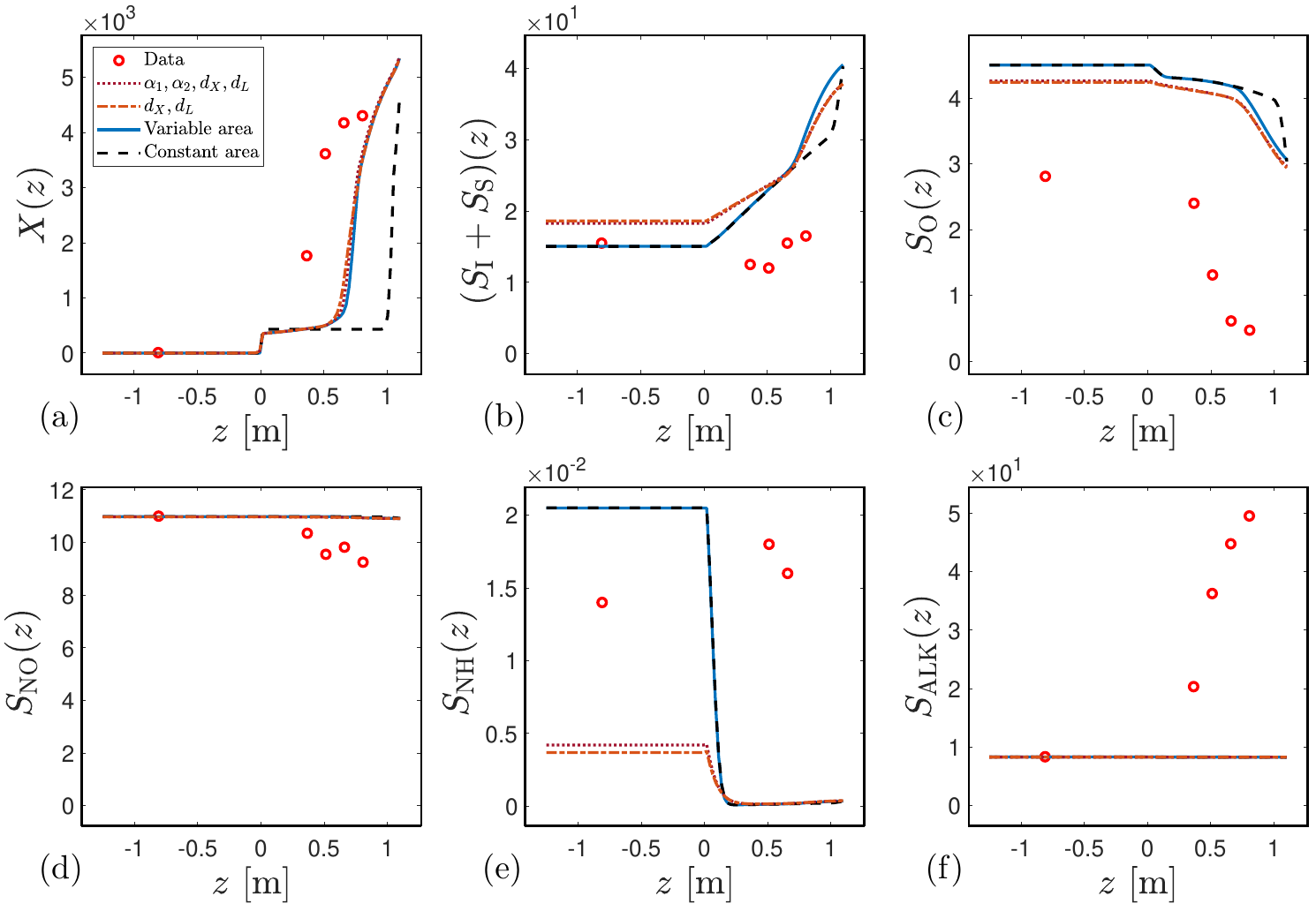} 
 \caption{Scenario L: Concentrations with units given in Table~\ref{table:AMS1_vari}. Each plot show five points of experimental data (circles) and simulation results at steady state with for the four variants of fitted model based only on Scenario~M.}  \label{fig:Scen1c}
\end{figure*}

\begin{figure*}[t!] 
\centering 
\includegraphics{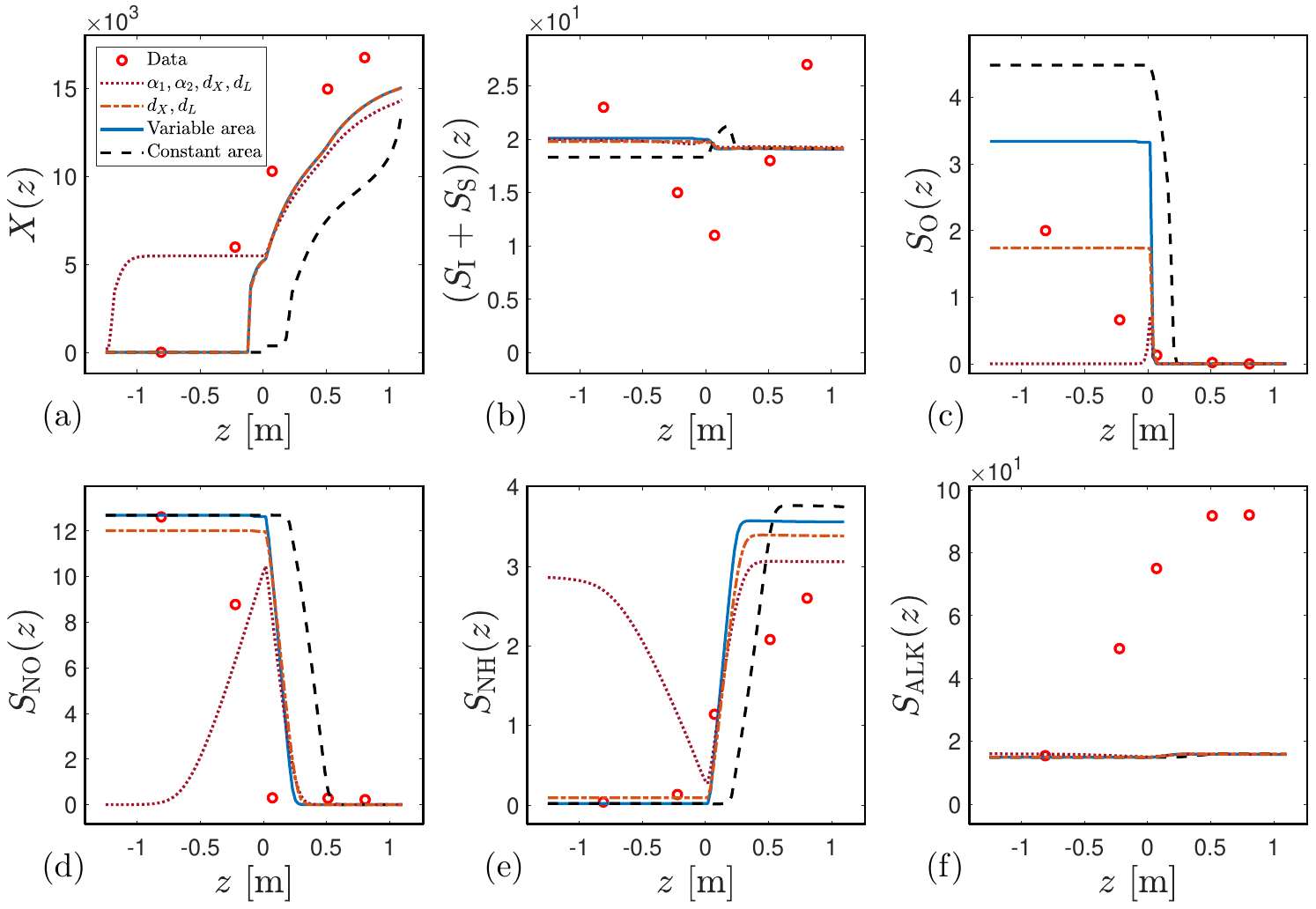} 
\caption{Scenario H: Concentrations with units given in Table~\ref{table:AMS1_vari}. Each plot show five points of experimental data (circles) and simulation results at steady state with for the four variants of fitted model based only on Scenario~M.} \label{fig:Scen3}
\end{figure*}

\begin{figure*}[t!] 
\centering 
\includegraphics{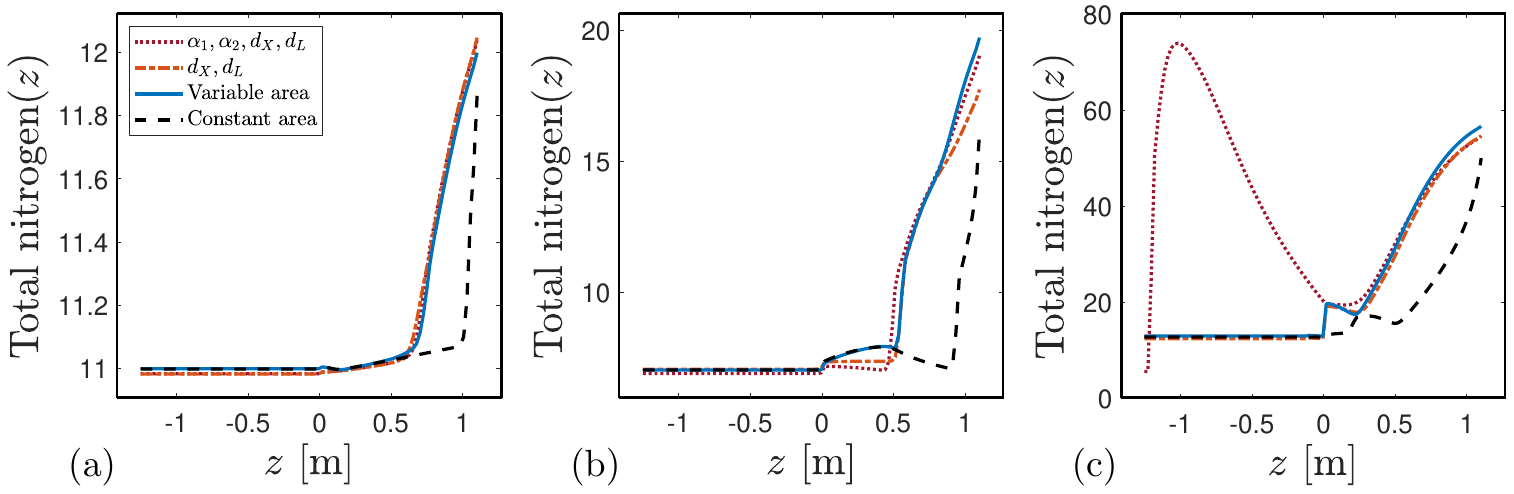} 
 \caption{Simulated steady states of total nitrogen for (a) Scenario~L, (b) Scenario~M and (c) Scenario~H with the four variants of fitted model.} \label{fig:Total_Nitrog}
\end{figure*}

As a validation, we keep the obtained parameter vales for Scenario~M and with these simulate Scenarios~L and H; see Figures~\ref{fig:Scen1c} and \ref{fig:Scen3}.
In addition, Figure~\ref{fig:Total_Nitrog} shows the simulated total nitrogen in the tank, which is the sum of nitrate, nitrite, ammonia and organic nitrogen, i.e., $X_{\rm ND}+S_{\rm NO}+S_{\rm NH}+S_{\rm ND}$.

The results indicate that a significant increase of simulation accuracy is achieved if the variability of the cross-sectional area is taken into account, even if the resulting model is still a spatially one-dimensional one.
The additional inclusion of the parameters for dispersion only $(d_X,d_L)$ improved the estimation of the concentration~$S_\mathrm{O}$ in all three scenarios, but the SBL only to a small degree.
The additional inclusion of the mixing parameters~$(\alpha_1,\alpha_2)$ caused only a small improvement in Scenario~M (as it should), whereas it did not improve the model prediction of Scenario~L, and in Scenario~H the prediction even got worse.
As Figure~\ref{fig:Scen3}~(a) shows, the simulated mixing near the inlet predicted too much sludge above the feed lever, and a bad prediction of the concentrations~$S_\mathrm{NO}$ in Figure~\ref{fig:Scen3}~(d) and $S_\mathrm{NH}$ in Figure~\ref{fig:Scen3}~(e).

\section{Conclusions}

A reactive-settling PDE model with standard parameters for a modified ASM1 model and several constitutive assumptions on the movement of particles and dissolved substrates have been calibrated to experimental data: 22~conventional batch sedimentation experiments and one steady-state scenario of an SST in a pilot plant.
The predictability of the model was evaluated to two further experimental SST scenarios.
The properties of hindered settling and compression at high concentrations for the flocculated particles was successfully fitted to 22~conventional batch sedimentation experiments; see Figure~\ref{fig:SBH}.
This was possible after the effects of the initial induction period of each test had been  transformed away.
One experimental steady-state scenario (four data points along the depth of the tank for each concentration) was thereafter used for the additional calibration of terms in the equations modelling hydrodynamic dispersion of partly the particles, and partly the dissolved substrates.
Adding more terms and parameters of course always leads to a better fit to the data used in the calibration.
Therefore, two experimental scenarios were used for validation.
Including only hydrodynamic dispersion lead to some improvement in the predictability of the model.
The additional inclusion of a term modelling the mixing of the suspension near the feed inlet lead to a worse predictability.
In contrast to the other phenomena, which were included in the derivation of the model, the inclusion of a general mixing term was made afterwards in an ad hoc way.

\section*{Acknowledgements}
Gamze Kirim, Elena Torfs and Peter Vanrolleghem are gratefully acknowledged for sharing their experimental data.
RB~and JC are   supported by ANID (Chile) through 
  Centro de Modelamiento Matem\'{a}tico  (CMM; BASAL pro\-ject FB210005). In addition 
 RB~is   supported by Fon\-de\-cyt project 1210610; 
  Anillo project ANID/PIA/ACT210030;  and CRHIAM, project   ANID/FON\-DAP/15130015. SD~ack\-nowledges support from the Swedish Research Council (Vetenskapsr\aa det, 2019-04601). RP  is supported by  scholarship ANID-PCHA/Doctorado Nacional/2020-21200939.

\bibliographystyle{elsarticle-harv}
\biboptions{authoryear}
\bibliography{ref_copy}

\appendix

\section{The modified ASM1 model}\label{sec:appendixA}
 \begin{table*}[t]
\caption{Stoichiometric and kinetic parameters.}
\label{table:Name_para}
\begin{center} 
\begin{tabular}{lp{10.5cm}ll} \hline 
 Symbol            & Name & Value & Unit \\ 
\hline
 $Y_{\rm A}$   & Yield for autotrophic biomass & 0.24 & $\rm (g\,COD)(g\,N)^{-1}$ \\
 $Y_{\rm H}$   & Yield for heterotrophic biomass & 0.57 & $\rm  (g\,COD)(g\,COD)^{-1}$ \\
 $f_{\rm P}$   & Fraction of biomass leading to particulate products & 0.1 & dimensionless\\
 $i_{\rm XB}$  & Mass of nitrogen per mass of COD in biomass & 0.07 & $\rm (g\,N)(g\,COD)^{-1}$  \\
 $i_{\rm XP}$  & Mass of nitrogen per mass of COD in products from biomass & 0.06 & $\rm (g\,N)(g\,COD)^{-1}$ \\
 $\mu_{\rm H}$ & Maximum specific growth rate for heterotrophic biomass & 4.0 & $\rm d^{-1}$ \\
 $K_{\rm S}$   & Half-saturation coefficient for heterotrophic biomass & 20.0 & $\rm (g\,COD)\,m^{-3}$ \\
 $K_{\rm O, H}$& Oxygen half-saturation coefficient for heterotrophic biomass & 0.25 & $\rm -(g\,COD)\,m^{-3}$ \\
 $K_{\rm NO}$  & Nitrate half-saturation coefficient for denitrifying heterotrophic biomass & 0.5 & $\rm (g\,NO_{3}\text{-}N)\, m^{-3}$ \\
 $b_{\rm H}$   & Decay coefficient for heterotrophic biomass & 0.5 & $\rm d^{-1}$ \\
 $\eta_{\rm g}$& Correction factor for $\mu_{\rm H}$ under anoxic conditions & 0.8 & dimensionless \\
 $\eta_{\rm h}$& Correction factor for hydrolysis under anoxic conditions & 0.35 & dimensionless \\
 $k_{\rm h}$   & Maximum specific hydrolysis rate & 1.5 & $\rm (g\,\text{COD})\, (g\,\text{COD})^{-1}\rm d^{-1}$ \\
 $K_{\rm X}$   & Half-saturation coefficient for hydrolysis of slowly biodegradable substrate & 0.02 & $\rm (g\,\text{COD})(g\,\text{COD})^{-1}$ \\ 
 $\mu_{\rm A}$ & Maximum specific growth rate for autotrophic biomass & 0.879 & $\rm d^{-1}$ \\
 $\bar{K}_{\rm NH}$ & Ammonia half-saturation coefficient for aerobic and anaerobic growth of heterotrophs & 0.007 & $\rm (g\,NH_{3}\text{-}N)\, m^{-3}$ \\
 $K_{\rm NH}$ & Ammonia half-saturation coefficient for autotrophic biomass & 1.0 &  $\rm (g\,NH_{3}\text{-}N)\, m^{-3}$\\
 $b_{\rm A}$   & Decay coefficient for autotrophic biomass & 0.132 & $\rm d^{-1}$ \\
 $K_{\rm O, A}$& Oxygen half-saturation coefficient for autotrophic biomass & 0.5 & $\rm -(g\,COD)\,m^{-3}$\\
 $k_{\rm a}$   & Ammonification rate & 0.08 & $\rm m^{3}(g COD)^{-1}d^{-1}$  \\
\hline 
\end{tabular} 
\end{center} 
\end{table*}

In terms of the stoichiometric matrices~$\bsigmaC$ and~$\bsigmaS$ and the vector~$\br(\boldsymbol{{C}},\boldsymbol{{S}})$ of eight processes of biokinetic reactions, the reaction rate vectors of \eqref{eq:modelC}, \eqref{eq:modelS} are 
\begin{align*} 
\bR_{\bC}(\boldsymbol{{C}},\boldsymbol{{S}}) & =\bsigmaC\br(\boldsymbol{{C}},\boldsymbol{{S}}),\\ 
\bR_{\bS}(\boldsymbol{{C}},\boldsymbol{{S}}) & =\bsigmaS\br(\boldsymbol{{C}},\boldsymbol{{S}}). 
 \end{align*} 
With the constants given in Table~\ref{table:Name_para}, the matrices are
{\footnotesize 
\begin{align*} 
 & \bsigmaC    \coloneqq   
\setlength\arraycolsep{2pt}
 \begin{bmatrix}
0&0&0&0&0&0&0&0 \\
0&0&0& 1-f_{\rm P} & 1-f_{\rm P} & 0 & -1 & 0 \\
1 & 1 & 0 & -1 & 0 & 0 & 0 & 0 \\
0 & 0 & 1 & 0 & -1 & 0&0&0 \\
0 & 0 & 0 & f_{\rm P} & f_{\rm P} & 0 & 0 & 0 \\
0 & 0 & 0 & i_{\rm XB}-f_{\rm P}i_{\rm XP} & i_{\rm XB}-f_{\rm P}i_{\rm XP} & 0 & 0 & -1
                  \end{bmatrix},    \\     
& \bsigmaS   \coloneqq    
\setlength\arraycolsep{2pt}
\begin{bmatrix}
0&0&0&0&0&0&0&0 \\[1.5mm] 
-\dfrac{1}{Y_{\rm H}} & -\dfrac{1}{Y_{\rm H}} &0&0&0&0& 1 & 0\\[3mm]
-\dfrac{1-Y_{\rm H}}{Y_{\rm H}} & 0 & -\dfrac{4.57-Y_{\rm A}}{Y_{\rm A}} &0&0&0&0&0 \\[3mm]
0& -\dfrac{1-Y_{\rm H}}{2.86Y_{\rm H}} & \dfrac{1}{Y_{\rm A}} &0&0&0&0&0 \\[3mm]
-i_{\rm XB} & -i_{\rm XB} & -i_{\rm XB}-\dfrac{1}{Y_{\rm A}} &0&0& 1 & 0 & 0 \\[1.5mm]
0&0&0&0&0& -1 & 0 & 1  \\[1.5mm]
-\dfrac{i_{\rm XB}}{14} & \dfrac{1-Y_{\rm H}}{40.04 Y_{\rm H}}-\dfrac{i_{\rm XB}}{14} & -\dfrac{i_{\rm XB}}{14}-\dfrac{1}{7 Y_{\rm A}} &0&0& \dfrac{1}{14} &0&0
\end{bmatrix}.
\end{align*}}
To describe~$\br(\boldsymbol{{C}},\boldsymbol{{S}})$, we define the functions 
\begin{align*}
&\mu_7(X_{\rm S},X_{\rm B, H}) \\ 
&\coloneqq \begin{cases}
0&\text{if $X_{\rm S}=0$ and $X_{\rm B, H}=0$,}\\
 \dfrac{X_{\rm S}X_{\rm B, H}}{K_{\rm X} X_{\rm B,H}+X_{\rm S}}&\text{otherwise,}
\end{cases}
\\
&\mu_8(X_{\rm B, H},X_{\rm N D}) \\ 
&\coloneqq \begin{cases}
0 &\text{if $X_{\rm S}=0$ and $X_{\rm B, H}=0$,}\\
 \dfrac{X_{\rm B, H}X_{\rm N D}}{K_{\rm X}X_{\rm B,H} + X_{\rm S}} &\text{otherwise.}
\end{cases}
\end{align*}
These functions are introduced to obtain well-defined expressions if any concentration is zero.
In the first two components of the vector~$\br$, we introduce an extra Monod factor with a small half-saturation parameter~$\smash{\bar{K}_{\rm NH}}$ for the concentration~$\smash{S_{\rm NH}}$ to guarantee that no consumption of ammonia can occur if its concentration is zero.
The rate vector is

\vspace{-0.5em}
{\footnotesize 
\begin{multline*} 
\vspace*{-10pt} \br(\boldsymbol{{C}},\boldsymbol{{S}})\coloneqq  \\ \vspace*{-10pt} 
\begin{pmatrix}
\mu_{\rm H} \dfrac{S_{\rm N H}}{\bar{K}_{\rm N H}+S_{\rm N H}} \dfrac{S_{\rm S}}{K_{\rm S}+S_{\rm S}} \dfrac{S_{\rm O}}{K_{\rm O, H}+S_{\rm O}} X_{\rm B, H}\\[4mm]
\mu_{\rm H} \dfrac{S_{\rm N H}}{\bar{K}_{\rm N H}+S_{\rm N H}} \dfrac{S_{\rm S}}{K_{\rm S}+S_{\rm S}} \dfrac{K_{\rm O, H}}{K_{\rm O, H}+S_{\rm O}} \dfrac{S_{\rm N O}}{K_{\rm N O}+S_{\rm N O}} \eta_{\rm g} X_{\rm B, H}\\[4mm]
\mu_{\rm A} \dfrac{S_{\rm N H}}{K_{\rm N H}+S_{\rm N H}} \dfrac{S_{\rm O}}{K_{\rm O, A}+S_{\rm O}} X_{\rm B, A}\\[3mm]
b_{\rm H} X_{\rm B, H}\\
b_{\rm A} X_{\rm B, A}\\
k_{\rm a} S_{\rm N D} X_{\rm B, H}\\[3mm]
k_{\rm h} \mu_7(X_{\rm S},X_{\rm B, H})\biggl( \dfrac{S_{\rm O}}{K_{\rm O, H}+S_{\rm O}}+\eta_{\rm h} \dfrac{K_{\rm O, H}}{K_{\rm O, H}+S_{\rm O}} \dfrac{S_{\rm N O}}{K_{\rm N O}+S_{\rm N O}}\biggr)\\[4mm] 
k_{\rm h} \mu_8(X_{\rm B, H},X_{\rm N D})\biggl(  \dfrac{S_{\rm O}}{K_{\rm O, H}+S_{\rm O}}+\eta_{\rm h} \dfrac{K_{\rm O, H}}{K_{\rm O, H}+S_{\rm O}} \dfrac{S_{\rm N O}}{K_{\rm N O}+S_{\rm N O}}\biggr)
\end{pmatrix}.
\end{multline*}}


\section{Numerical method}\label{sec:appendixB}

Combining ingredients from~\citep{SDIMA_MOL,bcdp_part2}, we suggest the following numerical method for the approximate solution of~\eqref{eq:modelC}, \eqref{eq:modelS}.

The height of the SST is divided into $N$~internal computational cells, or layers, of  depth $\Delta z = (B+H)/N$. 
The midpoint of layer~$j$ (numbered from above) has the coordinate $z=z_j$; hence, the layer is the interval $[z_{j-1/2},z_{j+1/2}]$ and we denote its average concentration vector by~$\smash{\bC_{j}(t)}$, which thus approximates $\bC(z_j,t)$, and similarly for $\boldsymbol{{S}}$ of the system~\eqref{eq:modelC}, \eqref{eq:modelS}
(we thus skip the tildes over numerical variables).
Recall that $X$ is always given by~\eqref{eq:X}.
The feed inlet at $z=0$ is located in the `feed layer'~$j_{\mathrm{f}}: = \lceil
H/\Delta z \rceil$, which is equal to the smallest integer larger than or equal to $H/\Delta z$. 
Above the interval $(-H,B)$, we add one layer to obtain the correct effluent concentrations via 
$\smash{\bC_\mathrm{e}(t)\coloneqq \bC_{0}(t)}$, 
and one layer below for the underflow concentration 
$\smash{\bC_\mathrm{u}(t)\coloneqq \bC_{N+1}(t)}$ 
 (analogously for~$\bS$).
For technical reasons, we set 
$\bC_{-1}\coloneqq  \bzero$  and $\bC_{N+2}\coloneqq \bzero$,  
and analogously for other variables. 
The cross-sectional area is approximated by
\begin{align*}
 A_{j+1/2}  &\coloneqq  \dfrac{1}{\Delta z}\int_{z_{j-1}}^{z_j} A(\xi)\,{\rm d} \xi 
 \end{align*} 
  and 
  \begin{align*}  A_j & \coloneqq  \dfrac{1}{\Delta z}\int_{z_{j-1/2}}^{z_{j+1/2}} A(\xi)\,{\rm d }\xi.
\end{align*}
We let $\gamma_{j+1/2}\coloneqq \gamma(z_{j+1/2})$ (similarly for other variables) and define the approximate volume average velocity
\begin{align*}
 q_{j+1/2} &\coloneqq \begin{cases}
-\Qe(t)/A_{j+1/2} & \text{for $j<j_{\mathrm{f}}$,} \\
\Qu(t)/A_{j+1/2} & \text{for $j\geq j_{\mathrm{f}}$.} 
\end{cases}
\end{align*}

With $\smash{\delta_{j,j_{\mathrm{f}}}}$ denoting the Kronecker delta, which is~1 if $j=j_{\mathrm{f}}$ and zero otherwise, the method-of-lines (MOL) formulation of the numerical method is
\begin{equation} \label{eq:MOL} \tag{B.1} 
\begin{aligned}
\frac{\rmd\bC_j}{\rmd t} &  = -\frac{\bPhi^{\bC}_{j+1/2}-\bPhi^{\bC}_{j-1/2}}{A_j\Delta z} +\delta_{j,j_{\mathrm{f}}}\frac{\bC_{\rm f}\Qf}{A_j\Delta z}
                           +\gamma_j\bR_{C,j},\\
\frac{\rmd\bS_j}{\rmd t} & = -\frac{\bPhi^{\bS}_{j+1/2}-\bPhi^{\bS}_{j-1/2}}{A_j\Delta z} +\delta_{j,j_{\mathrm{f}}}\frac{\bS_{\rm f}\Qf}{A_j\Delta z} 
                           + \gamma_j\bR_{S,j}. 
\end{aligned}
\end{equation}
Defining  $\smash{\chi(X)\coloneqq  \chi_{\{0<X <X_{\rm c} \}} (X)}$, 
$a^+ \coloneqq  \max \{0, a\}$ and $a^-\coloneqq \min \{ 0, a\}$,  
we may specify the numerical fluxes as follows. Utilizing the quantities  
\begin{align*} 
J_{j+1/2}^{\bC} &\coloneqq  \big(D(X_{j+1})-D(X_j)\big)/\Delta z,\\
J_{j+1/2}^{\rm d} &\coloneqq  \bigl(\chi (X_{j+1}) |q_{j+1/2}|\log(X_{j+1}) \\
& \qquad  - \chi (X_{j}) |q_{j-1/2}|\log(X_j)\bigr)/\Delta z,\\
v^{X}_{j+1/2}&\coloneqq  q_{j+1/2} + \gamma_{j+1/2}\big(\vhs(X_{j+1}) - J_{j+1/2}^{\bC}\big),\\
v^{L}_{j+1/2}& \coloneqq   \gamma_{j+1/2}(\vhs(X_{j+1}) - J_{j+1/2}^{\bC} - d_XJ_{j+1/2}^{\rm d}),\\
F^X_{j+1/2}&\coloneqq  v^{X,-}_{j+1/2} X_{j+1} + v^{X,+}_{j+1/2}  X_{j},\\
F^L_{j+1/2}&\coloneqq  v^{L,-}_{j+1/2} X_{j+1} + v^{L,+}_{j+1/2} X_{j},
\end{align*} 
we compute the numerical fluxes 
\begin{multline*} 
 \bPhi^{\bC}_{j+1/2}
\coloneqq  A_{j+1/2} \biggl(v^{X,-}_{j+1/2}\bC_{j+1} + v^{X,+}_{j+1/2} \bC_{j}\\
 -\gamma_{j+1/2}\bigl(\chi (X_{j+1}) |q_{j+1/2}| + d_{\mathrm{mix},j+1/2}\bigr) 
 \frac{\bC_{j+1}-\bC_{j}}{\Delta z}\biggr)
\end{multline*} 
and 
\begin{multline*} 
 \bPhi^{\bS}_{j+1/2} \\
 \coloneqq  A_{j+1/2}  \Bigg(\dfrac{((\rho_X-X_{j+1})q_{j+1/2}-F^L_{j+1/2})^-}{\rho_X-X_{j+1}}\bS_{j+1} \\  + \dfrac{((\rho_X-X_{j})q_{j+1/2}-F^L_{j+1/2})^+}{\rho_X-X_{j}}\bS_{j}\\
-\gamma_{j+1/2} \bigl(d_L|v^L_{j+1/2}| + d_{\mathrm{mix},j+1/2} \bigr)  \frac{\bS_{j+1}-\bS_{j}}{\Delta z}\Bigg).
\end{multline*}

Although any ODE solver can be used for the MOL system~\eqref{eq:MOL}, it is not meaningful  to use any higher order time-stepping algorithm since the spatial discretization is at most first-order accurate.
If $T$ is the simulation time, we let $t_n$, $n=0,1,\ldots,N_T$, denote the discrete time points and $\Delta t\coloneqq T/N_T$ the time step.
For explicit schemes, the right-hand sides of Equations~\eqref{eq:MOL} are evaluated at time~$t_n$.
The value of a variable at~$t_n$ is denoted by~$\smash{\bC_j^n}$, etc.\ and we set
\begin{align*}
Q_\mathrm{f}^n\coloneqq \frac{1}{\Delta t}\int_{t_n}^{t_{n+1}}Q_\mathrm{f}(t)\,\rmd t
\end{align*}
and similarly for the time-dependent reaction terms.
The time derivatives in \eqref{eq:MOL} are approximated by 
\begin{equation*}
\dd{\bC_j}{t}(t_n) \approx\frac{\bC_j^{n+1}-\bC_j^{n}}{\Delta t}.
\end{equation*}
This yields  the explicit scheme 
\begin{align*}
\bC_j^{n+1} & = \bC_j^n + \dfrac{\Delta t}{A_j\Delta z}
\bigl(-\bPhi^{\bC,n}_{j+1/2}+\bPhi^{\bC,n}_{j-1/2} \\ & \qquad \qquad  \qquad + \delta_{j,j_{\mathrm{f}}}\bC_{\rm f}^n Q_{\rm f}^n+\gamma_jA_j\Delta z\,\bR_{\bC,j}^n\bigr), \\
\bS_j^{n+1} & = \bS_j^n+\dfrac{\Delta t}{A_j\Delta z}
\bigl(-\bPhi^{\bS,n}_{j+1/2}+\bPhi^{\bS,n}_{j-1/2} \\ & \qquad \qquad  \qquad + \delta_{j,j_{\mathrm{f}}}\bS_{\rm f}^n Q_{\rm f}^n+\gamma_jA_j\Delta z\,\bR_{\bS,j}^n\bigr).
\end{align*}

\end{document}